\documentclass[12pt]{amsart}

\usepackage{latexsym,amssymb,amsmath}
\usepackage[dvips]{graphicx}

\def\ff{{\mathcal F}}

\def\mm{{\mathcal M}}

\def\ss{{\mathcal S}}
\def\ttt{{\mathcal T}}

\def\ffi{\varphi}

\def\dst{\displaystyle}

\def\supp{{\mathrm{supp}\,}}

\def\C{{\mathbb{C}}}
\def\E{{\mathbb{E}}}

\def\R{{\mathbb{R}}}

\def\T{{\mathbb{T}}}

\def\Z{{\mathbb{Z}}}

\def\cvh{{\natural\,}}
\newcommand{\norm}[1]{{\left\|{#1}\right\|}}

\newcommand{\abs}[1]{{\left|{#1}\right|}}
\newcommand{\scal}[1]{{\left\langle{#1}\right\rangle}}

\newenvironment{notation}[1][]{\vskip1pt\noindent\rm\textit{Notation}\,:\ }{\rm\vskip1pt}
\newenvironment{remark}[1][]{\vskip1pt\noindent\rm\textit{Remark #1}\,:\ }{\rm\vskip1pt}

\newtheorem{lemma}{Lemma}[section]

\newtheorem{theorem}[lemma]{Theorem}
\newtheorem{corollary}[lemma]{Corollary}

\begin{document}

 \baselineskip=17pt 

\title[A characterization of Fourier transforms]{A characterization of Fourier transforms}

\author{Philippe Jaming}

\address{Universit\'e d'Orl\'eans\\
Facult\'e des Sciences\\ 
MAPMO - F\'ed\'eration Denis Poisson\\ BP 6759\\ F 45067 Orl\'eans Cedex 2\\
France}
\email{Philippe.Jaming@univ-orleans.fr}

\begin{abstract}
The aim of this paper is to show that, in various situations, the only continuous linear map that transforms a convolution product into a pointwise product is a Fourier transform.
We focus on the cyclic groups $\Z/nZ$, the integers $\Z$, the Torus $\T$ and the real line.
We also ask a related question for the twisted convolution.
\end{abstract}

\subjclass{42A;A3842A85;42B10;43A25}

\keywords{Fourier transform;convolution}


\thanks{The author wishes to thank S. Alesker for sending him the preprint \cite{AAAV2}
and O. Guedon for pointing to that paper.
This work was partially financed by the French ANR project {\sl AHPI} (ANR-07-BLAN-0247-01)}
\maketitle
\hfill {\it In memory of A. Hulanicki.}

\section{Introduction}

The aim of this paper is to characterize the Fourier transform by some of its properties.
Indeed, the Fourier transform is well known to change a translation into
a modulation (multiplication by a character) and vice-versa and to change a convolution into a pointwise product. Moreover, these are some of its main features and are fundamental properties in many of its applications. The aim of this paper is to show that the Fourier transform is, to some extend,
uniquely determined by some of these properties.

Before going on, let us introduce some notation. Let $G$ be a locally compact Abelian group with Haar measure $\nu$ and let $\hat G$ be the dual group. Operations on $G$ will be denoted
additively. Let us recall that the convolution on $G$ is defined 
for $f,g\in L^1(G)$ by
$$
f*g(x)=\int_G f(t)g(x-t)\,\mbox{d}\nu(t)
$$
(and $f*g\in L^1(G)$) while the Fourier transform is defined by
$$
\ff(f)(\gamma)=\hat f(\gamma)=\int_G f(t)\overline{\gamma(t)}\,\mbox{d}\nu(t).
$$
We will here mainly focus on the four following cases,
$G=\hat G=\Z/n\Z$, $G=\Z$ and $\hat G=\T$ and vice versa or $G=\hat G=\R$,
(our results will then easily extend to products to such groups.)

We will here focus on two types of results. The first ones concerns the characterization of 
the Fourier transform as being essentially the only continuous linear transform that changes a convolution product into a pointwise product. To our knowledge the first results in that direction appear in the work of Lukacs \cite{Lu1,Lu2}, pursued in \cite{Em}, and an essentially complete result appeared in \cite{Fi}
for all LCA groups, under the mild additional constraint that the transform has a reasonable kernel.
We will show here that this hypothesis can be lifted. Further, a striking result, recently proved by
Alekser, Artstein-Avidan, Milman \cite{AAAV1,AAAV2} is that, to some extend, continuity and
linearity may be removed as well. More precisely, let us denote by $\ss(\R^d)$ the Schwartz functions on $\R^d$ and by $\ss'(\R^d)$ the Schwartz (tempered) distributions.

\medskip

\noindent{\bf Theorem (Alekser, Artstein-Avidan, Milman)}\\
{\sl Let $T\,:\ss(\R^d)\to\ss(\R^d)$ be a mapping that extends to a mapping
$T\,:\ss'(\R^d)\to\ss'(\R^d)$ that is bijective and such that
\begin{enumerate}
\renewcommand{\theenumi}{\roman{enumi}}
\item for every $f\in\ss(\R^d)$ and $g\in\ss'(\R^d)$, $T(f*g)=T(f).T(g)$;

\item for every $f\in\ss(\R^d)$ and $g\in\ss'(\R^d)$, $T(f.g)=T(f)*T(g)$.
\end{enumerate}
Then there exists $B\in\mm_n(\R)$ with $\det B=1$ such that $T(f)=\ff(f)\circ B$.}

\medskip

Note that $T$ is neither assumed to be linear nor to be continuous. We will adapt the proof
of this theorem to obtain an analogue result on the cyclic group. This has the advantage to highlight
the main features which come into proof of this theorem. The main difference is that in this theorem,
we assume that $T$ sends smooth functions into smooth functions. In the case of the cyclic group,
we do not have such functions at hand and are therefore lead to assume some mild continuity;
{\it see} Theorem \ref{th:A} for a precise statement.

\medskip

A second set of results has its origin in the work of Cooper \cite{Co1,Co2}. Here one considers
the Fourier transform as an intertwining operator between two groups of transforms
acting on $L^p$-spaces. In order to state the precise result, let us define,
for $\alpha\in\R$ and $f$ a function on $\R$, $\tau_\alpha f(t)=f(t+\alpha)$.
Further, if $\ffi\,:\R\to\R$, let $M^{(\ffi)}_\alpha f(t)=e^{i\alpha\ffi(t)}f(t)$.
It is easy to see that $\ff\tau_\alpha=M^{(t)}_\alpha\ff$ and
$\ff M^{(-t)}_\alpha=\tau_\alpha\ff$ {\it i.e.} the Fourier transform intertwines
translations and modulations and vice versa. The converse is also true. More precisely:

\medskip

\noindent{\bf Theorem (Cooper)}\\
{\sl Let $T\,: L^2(\R)\to L^2(\R)$ be a continuous linear transformation such that
there exists two measurable functions $\ffi,\psi\,:\R\to\R$ for which
$$
T\tau_\alpha=M^{(\ffi)}_\alpha T\quad\mathrm{and}\quad TM^{(\psi)}_\alpha=\tau_\alpha T.
$$
Then $\ffi(t)=bt+c$, $\psi(t)=bt+d$ with $b,c,d\in\R$ and $T=\ff$.}

\medskip

We will extend this theorem to $\Z/n\Z$ and $\Z$.

\medskip

The remaining of the article is organized as follows. In the next section, we will prove
the results for the groups $G=\Z$ and $G=\Z/n\Z$ while Section \ref{sec:3} is devoted to the cases
$G=\T^d$ and $G=\R^d$. We conclude with some questions concerning the twisted convolution.

\medskip

Before going on, let us introduce some more notation.
If $E\subset G$, we will denote by $\chi_E$ the function on $G$ given by $\chi_E(k)=1$ if $k\in E$ and $\chi_E(k)=0$ otherwise. The Kronecker symbol is denoted by $\delta_{j,k}$ where
$\delta_{j,k}=0$ or $1$ according to $j\not=k$ or $j=k$.

\section{The cyclic group and the integers}

In this section, we consider $G=\Z/n\Z$ or $G=\Z$.
We will write $\mathcal{C}(\hat G)$ for the set of $n$-periodic sequences
when $G=\hat G=\Z/n\Z$ or of continuous functions on $\hat G=\T$ if $G=\Z$.
Our first result is the following:

\begin{theorem}\label{th:1}
Let $G=\Z/n\Z$ or $G=\Z$.
Let $T$ be a linear continuous map $T\,:L^1(G)\to \mathcal{C}(\hat G)$ such that
$T(f*g)=T(f).T(g)$. Then there exists $E\subset \hat G$ and a map 
$\sigma\,:\hat G\to\hat G$ such that, for
$f\in L^1(G)$ and almost every $\eta\in\hat G$
$T(f)(\eta)=\chi_E(\eta)\widehat{f}\bigl(\sigma(\eta)\bigr)$. *Moreover, $\sigma$ is measurable if $G=\Z$.
\end{theorem}

\begin{proof}
Let $\delta_k=(\delta_{j,k})_{j \in G}\in L^1(G)$. Then $\delta_k*\delta_\ell=\delta_{k+\ell}$, so that
\begin{equation}
\label{eq:Tdelta}
T(\delta_{k+\ell})=T(\delta_k*\delta_\ell)=T(\delta_k)T(\delta_\ell).
\end{equation}
In particular, for each $\eta\in \hat G$, the map $\pi_\eta\,:k\to T(\delta_k)(\eta)$ is a group homomorphism from $G$ to $\C$.

First note $\pi_\eta(0)=\pi_\eta(k)\pi_\eta(-k)$
so that if $\pi_\eta$ vanishes somewhere, it vanishes at $0$. Conversely 
$\pi_\eta(k)=\pi_\eta(k)\pi_\eta(0)$
so that if $\pi_\eta$ vanishes at $0$, it vanishes everywhere. Further 
$\pi_\eta(0)=\pi_\eta(0)^2$ so that $\pi_\eta(0)=0$ or $1$. 

We will now assume that $\pi_\eta(0)=1$ and exploit 
$\pi_\eta(k+1)=\pi_\eta(k)\pi_\eta(1)$ which implies that
$\pi_\eta(k)=\pi_\eta(1)^k$.
We now need to distinguish the two cases:

\smallskip

--- if $G=\Z/n\Z$, then $1=\pi_\eta(0)=\pi_\eta(n)=\pi_\eta(1)^n$
$\pi_\eta(1)$ is an $n$-th root of unity {\it i.e.}
$T(\delta_1)(\eta)=e^{2i\pi \sigma(\eta)/n}$ for some $\dst\sigma(\eta)\in \{0,1\ldots,n-1\}=\Z/n\Z$. 
In follows that $T(\delta_k)(\eta)=e^{2i\pi k\sigma(\eta)/n}$. 

\smallskip

--- if $G=\Z$, as  $T$ was assumed to be continuous $L^1(G)\to \mathcal{G}(\hat G)$, there is a constant $C>0$ such that, for every $f\in L^1(G)$, $\norm{Tf}_\infty\leq C\norm{f}_1$.
In particular, for every $k\in\Z$ and every $m\in\hat G=\T$
$$
|\pi_\eta(1)^k|=|[T(\delta_1)(\eta)]^k|=
|T(\delta_k)(\eta)|\leq C\norm{\delta_k}_1=C
$$
thus, by letting $k\to\pm\infty$, we obtain that $\pi_\eta(1)$ is a complex number of modulus $1$
(as it is not $0$ since $\pi_\eta(0)^k\not=0$). We may thus write
$T(\delta_1)(\eta)=e^{2i\pi \sigma(\eta)}$ for some $\sigma(\eta)\in[0,1]\simeq\T=\hat G$.
Moreover, as $\eta\to T(\delta_1)(\eta)$ is measurable, we may assume that $\sigma$ is measurable as well.

\smallskip

Let us now define $E=\{\eta\in \hat G\,:\ T(\delta_k)(\eta)=0\ \forall\ k\in G\}$. Then, by linearity and 
continuity of $T$, for $f\in L^1(G)$,
\begin{eqnarray*}
Tf(\eta)&=&T\left(\sum_{k\in G}f(k)\delta_k\right)(\eta)=\sum_{k\in G}f(k)T(\delta_k)(\eta)\\
&=&\begin{cases}\dst\sum_{k\in G}f(k)\chi_E(\eta)e^{2i\pi k\sigma(\eta)/n}&\mbox{if }G=\Z/n\Z\\
\dst\sum_{k\in G}f(k)\chi_E(\eta)e^{2i\pi k\sigma(\eta)}&\mbox{if }G=\Z\end{cases}\\
&=&\chi_E(\eta)\widehat{f}\bigl(\sigma(\eta)\bigr),
\end{eqnarray*}
which completes the proof.
\end{proof}

\begin{remark}
Using tensorization, we may extend the result with no difficulty to $G=\prod_{i\in I}\Z/n_i\Z\times\Z^d$.
\end{remark}

We will now adapt the proof of \cite{AAAV1,AAAV2} to show that on $\Z_n$, a bijective transform
that maps a product into a convolution is essentially a Fourier transform.
We will need some notation:

\begin{notation}
We will consider the following particular
elements of $L^1(\Z/n\Z)$\,: $\mathbf{0}=(0,\ldots,0)$ and $\mathbf{1}=(1,\ldots,1)$.
Further, if $a\in L^1(\Z/n\Z)$ we will write $\dst\E[a]=\sum_{j=0}^{n-1}a_j$.
\end{notation}

We can now state the main theorem:

\begin{theorem}\label{th:A}
Let $\ttt\,:L^1(\Z/n\Z)\to L^1(\Z/n\Z)$ be a bijective transformation (not necessarily linear) such
that the map $\C\to L^1(\Z/n\Z)$, $c\to \ttt(c\mathbf{1})$ is continuous. Assume that
\begin{enumerate}
\renewcommand{\theenumi}{\roman{enumi}}
\item for every $a,b\in L^1(\Z/n\Z)$, $\ttt(a.b)=\ttt(a).\ttt(b)$;
\item for every $a,b\in L^1(\Z/n\Z)$, $\ttt(a*b)=\ttt(a)*\ttt(b)$.
\end{enumerate}
Then there exists $\eta\in\{1,\ldots,n-1\}$ that has no common divisor with $n$ such that

-- either, for every $j\in\Z/n\Z$ and every $a\in L^1(\Z/n\Z)$, $\ttt(a)(\eta j)=a(j)$;

-- or, for every $j\in\Z/n\Z$ and every $a\in L^1(\Z/n\Z)$, $\ttt(a)(\eta j)=\overline{a(j)}$.
\end{theorem}

\begin{remark}
The fact that $\eta$ has no common divisor with $n$ implies that the map $j\to j\eta$
is a permutation of $\{0,\ldots,n-1\}$ so that the map $\ttt$ is actually fully determined.
\end{remark}

\begin{corollary}\label{cor:A}
Let $\ttt\,:L^1(\Z/n\Z)\to L^1(\Z/n\Z)$ be a bijective transformation (not necessarily linear) such
that the map $\C\to L^1(\Z/n\Z)$, $c\to T(c\mathbf{1})$ is continuous. Assume that
\begin{enumerate}
\renewcommand{\theenumi}{\roman{enumi}}
\item for every $a,b\in L^1(\Z/n\Z)$, $\ttt(a.b)=\ttt(a)*\ttt(b)$;
\item for every $a,b\in L^1(\Z/n\Z)$, $\ttt(a*b)=\ttt(a).\ttt(b)$.
\end{enumerate}
Then there exists $\eta\in\{1,\ldots,n-1\}$ that has no common divisor with $n$ such that
either, for every $j\in\Z/n\Z$ and every $a\in L^1(\Z/n\Z)$, $\ttt(a)(\eta j)=\hat a(j)$
or, for every $j\in\Z/n\Z$ and every $a\in L^1(\Z/n\Z)$, $\ttt(a)(\eta j)=\overline{\hat a(j)}$.
\end{corollary}

\begin{proof}[Proof of Corollary \ref{cor:A}]
It is enough to apply Theorem \ref{th:A} to  $\tilde \ttt=\ff^{-1}\ttt$.
\end{proof}

\begin{corollary}\label{cor:B}
Let $\ttt\,:L^1(\Z/n\Z)\to L^1(\Z/n\Z)$ be a bijective transformation (not necessarily linear) such
that the map $\C\to L^1(\Z/n\Z)$, $c\to \ttt(c\mathbf{1})$ is continuous. Assume that for
every $a\in L^1(\Z/n\Z)$, $\ttt^2a(k)=a(-k)$ and that one of the following two identities holds:
\begin{enumerate}
\renewcommand{\theenumi}{\roman{enumi}}
\item for every $a,b\in L^1(\Z/n\Z)$, $\ttt(a.b)=\ttt(a)*\ttt(b)$;
\item for every $a,b\in L^1(\Z/n\Z)$, $\ttt(a*b)=\ttt(a).\ttt(b)$.
\end{enumerate}
Then there exists $\eta\in\{1,\ldots,n-1\}$ that has no common divisor with $n$ such that
either, for every $j\in\Z/n\Z$ and every $a\in L^1(\Z/n\Z)$, $\ttt(a)(\eta j)=\hat a(j)$
or, for every $j\in\Z/n\Z$ and every $a\in L^1(\Z/n\Z)$, $\ttt(a)(\eta j)=\overline{\hat a(j)}$.
\end{corollary}

\begin{proof}[Proof of Corollary \ref{cor:B}]
If $\ttt^2a(k)=a(-k)$ then if one of the identities holds, so does the second one, so that
Corollary \ref{cor:A} gives the result.
\end{proof}

\begin{proof}[Proof of Theorem \ref{th:A}]
The proof goes in several steps that are similar to those in \cite{AAAV2}. The first one consists in identifying
the image by $\ttt$ of some particular elements of $L^1(\Z/n\Z)$:

\smallskip

\noindent{\bf Step 1.} {\sl We have $\ttt(\delta_0)=\delta_0$, $\ttt(\mathbf{0})=\mathbf{0}$ and $\ttt(\mathbf{1})=\ttt(\mathbf{1})$.
Moreover, there is a $k\in\{-1,1\}$ and an $\alpha\in\C$ with $\mbox{Re}\,\alpha>0$ such that,
if we define $\beta\,:\C\to\C$ by $\beta(0)=0$ and $\dst\beta(c)=\left(\frac{c}{|c|}\right)^k|c|^\alpha$ for $c\not=0$, then $\ttt(c\mathbf{1})=\beta(c)\mathbf{1}$.}

\smallskip

Indeed, as $\ttt(a.b)=\ttt(a).\ttt(b)$, we immediately get the following:
$$
\ttt(c_1c_2\mathbf{1})=\ttt(c_1\mathbf{1}).\ttt(c_2\mathbf{1})\quad\mbox{and}\quad
\ttt(c_1\delta_j)=\ttt(c_1\mathbf{1}).\ttt(\delta_j)
$$
while from $\ttt(a*b)=\ttt(a)*\ttt(b)$ we deduce that
$$
\ttt(\delta_{j+k})=\ttt(\delta_j)\ttt(\delta_k)\quad\mbox{and}\quad
\ttt(a)=\ttt(\delta_0*a)=\ttt(\delta_0)*\ttt(a).
$$
Applying this last identity to $a=\ttt^{-1}(\delta_0)$ we get $\delta_0=\ttt(\delta_0)*\delta_0=\ttt(\delta_0)$. 

Further $a=a.\mathbf{1}$ thus $\ttt(a)=\ttt(a).\ttt(\mathbf{1})$
and, applying this again to $a=\ttt^{-1}(b)$, we have $b=b.\ttt(\mathbf{1})$ for all $b\in\ell^2_n$, thus
$\ttt(\mathbf{1})=\mathbf{1}$.
Similarly, $\mathbf{0}=a.\mathbf{0}$ thus $\ttt(\mathbf{0})=\ttt(a).\ttt(\mathbf{0})$ and, applying this to $a=\ttt^{-1}(\mathbf{0})$ we get $\ttt(\mathbf{0})=\mathbf{0}.\ttt(\mathbf{0})=\mathbf{0}$.

Finally, $\E[a]\mathbf{1}=a*1$ thus $\ttt(\E[a]\mathbf{1})=\ttt(a*\mathbf{1})=\ttt(a)*\mathbf{1}=\E[\ttt(a)]\mathbf{1}$.
As every $c\in\C$ may be written $c=E[\frac{c}{n}\mathbf{1}]$, we may define
$\beta(c)=E\left[\ttt\left(\frac{c}{n}\mathbf{1}\right)\right]$
so that $T(c\mathbf{1})=\beta(c)\mathbf{1}$. Note that $\beta$ is continuous
since we have assumed that $T$ acts continuously on constants and as $\ttt$ is one-to-one, so is $\beta$.
Moreover, $\beta$ is multiplicative\,:
$$
\beta(c_1c_2)\mathbf{1}=\ttt(c_1c_2\mathbf{1})
=\ttt(c_1\mathbf{1}).\ttt(c_2\mathbf{1})=\beta(c_1)\beta(c_2)\mathbf{1}.
$$
It is then easy to check that there is a $k\in\{-1,1\}$ and an $\alpha\in\C$ with $\mbox{Re}\,\alpha>0$
such that $\beta(0)=0$ and $\beta(c)=\dst\left(\frac{c}{|c|}\right)^k|c|^\alpha$.

We will now take care of the image of $\delta_j$, $j=0,\ldots,n-1$.

\smallskip

\noindent{\bf Step 2.} {\sl There is an $\eta\in\{1,\ldots,n-1\}$ with no common divisor with $n$
such that $\ttt(\delta_j)=\delta_{\eta j}$.}

\smallskip

%
Assume that $k\not=\ell\in\supp \ttt(\delta_j)$ thus $\delta_k.\ttt(\delta_j)\not=\mathbf{0}$
and $\delta_\ell.\ttt(\delta_j)\not=\mathbf{0}$. Let $a=\ttt^{-1}(\delta_k)$,
$b=\ttt^{-1}(\delta_\ell)$, then
$$
a.\delta_j=\ttt^{-1}(\delta_k).\ttt^{-1}\bigl(\ttt(\delta_j)\bigr)
=\ttt^{-1}\bigl(\delta_k.\ttt(\delta_j))\not=\ttt^{-1}(\mathbf{0})
$$
since $\ttt$ is one-to-one. From Step 1, we know that $\ttt^{-1}(\mathbf{0})=\mathbf{0}$,
therefore $a.\delta_j\not=\mathbf{0}$. 
For the same reason, $b.\delta_j\not=\mathbf{0}$. In particular,
$a.b\not=\mathbf{0}$, thus $\ttt(a.b)\not=\mathbf{0}$. But this contradicts $\ttt(a.b)=\delta_k.\delta_\ell$ with $k\not=\ell$.

It follows that, for each $j\in\{1,\ldots,n-1\}$, there exists $c_j\in\C\setminus\{0\}$
and $\sigma(j)\in\{0,\ldots,n-1\}$ such that $\ttt(\delta_j)=c_j\delta_{\sigma(j)}$.
But then
$$
\mathbf{1}=\ttt(\mathbf{1})=\ttt(\mathbf{1}*\delta_j)=\ttt(\mathbf{1})*\ttt\delta_j
=c_j\mathbf{1}*\delta_j=c_j\mathbf{1}
$$
thus $c_j=1$. As $\ttt$ is one-to-one, it follows that $\sigma(j)\in\{1,\ldots,n-1\}$
and that $\sigma$ is a permutation.

Next
$$
\delta_{\sigma(j+k)}=\ttt(\delta_{j+k})=\ttt(\delta_j*\delta_k)=\ttt(\delta_j)*\ttt(\delta_k)=\delta_{\sigma(j)}*
\delta_{\sigma(k)}=\delta_{\sigma(j)+\sigma(k)}.
$$
Thus $\sigma(j+k)=\sigma(j)+\sigma(k)$ and therefore $\sigma(j)=j\sigma(1)$. Further, the
fact that $\sigma$ is a permutation implies that $\sigma(1)$ has no common divisor with $n$
(Bezout's Theorem).

\smallskip

\noindent{\bf Step 3.} {\sl Conclusion.}

\smallskip

We can now prove that $\ttt$ is of the expected form: fix $j\in\{0,\ldots,n-1\}$ and
$a\in\ell^2_n$. Let $k=\sigma^{-1}(j)$ so that $\ttt(\delta_k)=\delta_j$. Then
\begin{eqnarray*}
\ttt(a)(j)\delta_j&=&\ttt(a).\delta_j=\ttt(a).\ttt(\delta_k)=\ttt(a.\delta_k)\\
&=&\ttt\bigl(a(k)\mathbf{1}.\delta_k\bigr)
=\beta\bigl(a(k)\bigr)\mathbf{1}.\ttt(\delta_k)\\
&=&\beta\bigl(a(k)\bigr)\delta_j.
\end{eqnarray*}
It follows that $\ttt(a)(j)=\beta\bigl(a\circ \sigma^{-1}(j)\bigr)=
\dst\left(\frac{a\circ \sigma^{-1}(j)}{\abs{a\circ \sigma^{-1}(j)}}\right)^k|a\circ \sigma^{-1}(j)|^\alpha$. We want to prove that $\alpha=1$. But
\begin{eqnarray*}
\E[\ttt(a)]\mathbf{1}&=&\ttt(a)*\mathbf{1}=\ttt(a)*\ttt(\mathbf{1})=\ttt(a*\mathbf{1})=\ttt(\E[a]\mathbf{1})\\
&=&\left(\frac{\E[a]}{\overline{\E[a]}}\right)^k|\E[a]|^\alpha\mathbf{1}
\end{eqnarray*}
so that
$\E[\ttt(a)]=\left(\frac{\E[a]}{\abs{\E[a]}}\right)^k|\E[a]|^\alpha$ or, in other words,
\begin{eqnarray*}
\sum_{\ell=0}^{n-1}\left(\frac{a(\ell)}{\abs{a(\ell)}}\right)^k\abs{a(\ell)}^\alpha
&=&\sum_{j=0}^{n-1}\left(\frac{a\bigl(\sigma^{-1}(j)\bigr)}{\abs{a\bigl(\sigma^{-1}(j)\bigr)}}\right)^k
\abs{a\bigl(\sigma^{-1}(j)\bigr)}^\alpha\\
&=&
\left(\frac{\sum_{j=0}^{n-1}a_j}{\sum_{j=0}^{n-1}\abs{a_j}}\right)^k\abs{\sum_{j=0}^{n-1}a_j}^\alpha.
\end{eqnarray*}
Taking $a(0)=1$, $a(1)=t>0$ and $a(j)=0$ for $j=2,\ldots,n-1$, this reduces to
$1+t^\alpha=(1+t)^\alpha$. This implies that $\alpha=1$ (which is easiest seen by differentiating and letting $t\to0$). It follows that $\beta(c)=c$ or $\bar c$ according to $k=1$ or $-1$.
\end{proof}

\begin{remark}
The proof adapts with no difficulty to any finite Abelian group. To prove the same result on $\Z$,
the best is to first compose $\ttt$ with a Fourier transform and then to adapt the proof
in \cite{AAAV2} from the real line to the torus. We refrain from giving the details here.

The lines of proof given here follows those given in \cite{AAAV2} (up to the ordering and the removal of technicalities that are useless in the finite group setting). The main difference is that we need to assume
that $\ttt$ acts continuously on constants. In \cite{AAAV2} this hypothesis is replaced by the fact that
$\ttt$ sends smooth functions into smooth functions.

Finally, it should also be noted that Hypothesis (i) and (ii) are only used when either $a$ or $b$
is either a constant $c\mathbf{1}$ or a Dirac $\delta_j$.
\end{remark}

We will conclude this section with a Cooper like theorem. Let us first introduce some notation:

\begin{notation}
Let $\ffi\,:\Z/n\Z\to \C$. For $k\in\Z/n\Z$, we define the two following linear operators $L^1(\Z/n\Z)\to L^1(\Z/n\Z)$\,:
$$
\tau_ka(j)=a(j+k)\quad\mbox{and}\quad M^{(\ffi)}_k a(j)=e^{k\ffi(j)}a(j).
$$
Note that actually $\ffi\,:\Z/n\Z\to \C/2i\pi\Z$.
As is well known, if $\ffi(j)=2i\pi j/n$ for some $k\in\Z/n\Z$, $\ff\tau_{-k}=M^{(\ffi)}_k\ff$ and
$\ff M^{(\ffi)}_k=\tau_k\ff$.
\end{notation}

We can now state the following:

\begin{theorem}
Let $\ttt\,:L^1(\Z/n\Z)\to L^1(\Z/n\Z)$ be continuous linear operator such that there exist
two maps $\ffi,\psi\,:\Z/n\Z\to \C$ for which
$$
\ttt\tau_k=M_k^{(\ffi)}\ttt\quad\mbox{and}\quad \ttt M_k^{(\psi)}=\tau_k\ttt.
$$
Then there exist $k_0,m_0,m_1\in\Z/n\Z$, $c\in\C$ such that
$\ffi(j)=\dst\frac{2i\pi}{n}(k_0j+m_0)$, $\psi(j)=\frac{2i\pi}{n}(-k_0j+m_1)$ and
$$
\ttt(a)(\ell)=ce^{2i\pi\ell m_1/n}\hat a(k_0\ell+m_0).
$$
\end{theorem}

\begin{proof} Without loss of generality, we may assume that $\ttt\not=0$.
First note that the conditions are equivalent to
\begin{equation}
\label{eq:cooper}
a)\ \ttt(\delta_k*a)(\ell)=e^{-k\ffi(\ell)}\ttt(a)(\ell)\quad\mbox{and}\quad b)\ \ttt(e^{-k\psi(\cdot)}a)=\delta_k*\ttt(a).
\end{equation}
Note that these two expressions are $n$-periodic in $k$ so that $\ffi$ and $\psi$ take their values
in $\{0,\frac{2i\pi}{n},\ldots,\frac{2i\pi(n-1)}{n}\}$.

First a) implies
$$
\ttt(\delta_j)(\ell)=\ttt(\delta_j*\delta_0)(\ell)=e^{-j\ffi(\ell)}\ttt(\delta_0)(\ell).
$$
Next \eqref{eq:cooper} b) implies that
\begin{eqnarray*}
e^{-k\psi(j)}\ttt(\delta_j)(\ell)&=&T(e^{-k\psi(j)}\delta_j)(\ell)=T(e^{-k\psi(\cdot)}\delta_j)(\ell)\\
&=&\delta_k*T(\delta_j)(\ell)=T(\delta_j)(\ell-k).
\end{eqnarray*}
In particular, $\ttt(\delta_j)(\ell)=e^{\ell\psi(j)}\ttt(\delta_j)(0)$, thus
$$
\ttt(\delta_j)(\ell)=e^{\ell\psi(j)-j\ffi(0)}\ttt(\delta_0)(0).
$$
From linearity, we thus get that for $a\in\ell^2_n$
$$
\ttt(a)(\ell)=\sum_{j=0}^{n-1}a(j)\ttt(\delta_j)(\ell)=
\left(\sum_{j=0}^{n-1}a(j)e^{\ell\psi(j)-j\ffi(0)}\right)\ttt(\delta_0)(0).
$$
As we assumed that $\ttt\not=0$, we thus have $\ttt(\delta_0)(0)\not=0$. Then \eqref{eq:cooper} reads
$$
\sum_{j=0}^{n-1}a(j)e^{\ell\psi(j+k)-(j+k)\ffi(0)}=
\sum_{j=0}^{n-1}a(j)e^{\ell\psi(j)-j\ffi(0)-k\ffi(\ell)}
$$
thus $\ell\psi(j+k)-(j+k)\ffi(0)=\ell\psi(j)-j\ffi(0)-k\ffi(\ell)$ for all $j,k,\ell\in\Z/n\Z$
(modulo $2i\pi/n$).
Taking $k=1$, we get
$$
\ffi(\ell)-\ffi(0)=\bigl(\psi(j)-\psi(j+1)\bigr)\ell
$$
so that $\ffi$ and $\psi$ are ``affine''. More precisely, $\ffi(\ell)=\bigl(\psi(0)-\psi(1)\bigr)\ell+\ffi(0)$ modulo $2i\pi/n$ and, as $\ffi$ takes its values
in $\frac{2i\pi}{n}\Z/n\Z$, 
$\ffi(\ell)=\frac{2i\pi}{n}(k_0\ell+m_0)$ (modulo $2i\pi/n$) with $k_0,m_0\in\{0,\ldots,n-1\}$ and $b\in\C$.
Further $\psi(j+1)=\psi(j)+\ffi(0)-\ffi(1)$ thus $\psi(j)=\psi(0)+j\bigl(\ffi(0)-\ffi(1)\bigr)=\frac{2i\pi }{n}(-k_0j+m_1)$ (again modulo $2i\pi/n$).

We thus conclude that
$$
\ttt(a)(\ell)=e^{2i\pi\frac{\ell m_1}{n}}\sum_{j=0}^{n-1}a(j)e^{-2i\pi\frac{k_0\ell+m_0}{n}j}.
$$
as expected.
\end{proof}

\section{The real line and the Torus}\label{sec:3}

We now consider the case $G=\R^d$ resp. $G=\T^d$ so that $\hat G=\R^d$ resp. $\hat G=\Z^d$. To simplify notation, we write $\mathcal{C}(\Z^d)=L^\infty(\Z^d)$.

\begin{theorem}\label{prop:fourR}
Let $d\geq 1$ be an integer and $G=\R^d$ or $G=\T^d$.
Let $T$ be a continuous linear operator $L^1(G)\to \mathcal{C}(\hat G)$ such that
$T(f*g)=T(f)T(g)$. Then there exists a set $E\subset G$ and a function
$\ffi\,: \hat G\to \hat G$ such that $T(f)(\xi)=\chi_E(\xi)\widehat{f}\bigr(\ffi(\xi)\bigl)$.
\end{theorem}

\begin{proof}
Let us fix $\xi\in\hat G$ and consider the continuous linear functional $T_\xi$
on $L^1(G)$ given by $T_\xi(f)=T(f)(\xi)$. Then there exists a bounded function $h_\xi$ on $G$
such that $T_\xi(f)=\int_G f(t)h_\xi(t)\,\mbox{d}t$. There is no loss of generality in assuming that
$h_\xi\not=0$.

Let us now take $A,B$ to sets of finite measure. Then Fubini's Theorem implies that
\begin{eqnarray}
\int_{A\times B}h_\xi(s+t)\,\mbox{d}s\,\mbox{d}t&=&\int_{\R^d} \chi_A*\chi_B(t)h_\xi(t)\,\mbox{d}t
=T(\chi_A*\chi_B)(\xi)\nonumber\\
&=&T(\chi_A)(\xi)T(\chi_B)(\xi)=\int_Ah_\xi(t)\,\mbox{d}t\,\int_Bh_\xi(t)\,\mbox{d}t.\label{eq:ffi}
\end{eqnarray}

Now let $\ffi_n$ be defined on $G^2$ by
$$
\ffi_n(x,y)=\begin{cases}\bigl(h_\xi(x+y)-h_\xi(x)h_\xi(y)\bigr)\chi_{[-n,n]}(x)\chi_{[-n,n]}(y)
&\mbox{if }G=\R\\
h_\xi(x+y)-h_\xi(x)h_\xi(y)&\mbox{if }G=\T\end{cases}.
$$
As $\ffi_n$ is bounded (since $h_\xi$ is) and has compact support, $\ffi_n\in L^1(G^2)$
and \eqref{eq:ffi} implies that
$$
\int_{A\times B}\ffi_n(x,y)\,\mbox{d}x\,\mbox{d}y=0
$$
for every sets $A,B$ of finite measure, so that $\ffi_n=0$ for every $n$. That is
\begin{equation}
\label{eq:h}
h_\xi(x+y)=h_\xi(x)h_\xi(y)\qquad\mbox{for almost every }x,y\in G.
\end{equation}
If $h_\xi$ were continuous, this would imply that $h_\xi(x)=e^{i\scal{a_\xi,x}}$ and, by boundedness of $h_\xi$,
that $a_\xi\in\R^d$. We will now overcome this difficulty by introducing
$$
H_{\xi,j}(x)=\int_0^xh_\xi(t\mathbf{e}_j)\,\mbox{d}t
$$
where $j=1,\ldots,d$ and $\mathbf{e}_j=(\delta_{j,k})_{k=1,\ldots,d}$ is the $j$-th vector in the standard basis.
Clearly $H_{\xi,j}$ is continuous and satisfies
$$
H_{\xi,j}(x)H_{\xi,j}(y)=\int_0^x\bigl(H_{\xi,j}(y+t)-H_{\xi,j}(t)\bigr)\,\mbox{d}t.
$$
From this, we immediately deduce that $H_{\xi,j}$ is smooth, that $H_{\xi,j}^\prime(t)=h_\xi(t\mathbf{e}_j)$
almost everywhere and that $H_{\xi,j}^\prime(x+y)=H_{\xi,j}^\prime(x)H_{\xi,j}^\prime(y)$ everywhere. Thus, for almost every $x\in\R$ or $\T$,
$h_\xi(x\mathbf{e}_j)=e^{ia_{\xi,j} x}$ with $a_{\xi,j}$ real. Finally for $x\in G$,
$$
h_\xi(x)=h_\xi(x_1\mathbf{e}_1+\cdots+x_d\mathbf{e}_d)=h_\xi(x_1\mathbf{e}_1)\cdots h_\xi(x_d\mathbf{e}_d)
=e^{i\scal{a_{\xi},x}}
$$
where $a_\xi=(a_{\xi,1},\ldots,a_{\xi,d})$.

We have thus proved that there exists a map from $\ffi\,:G\to G$ and a set $E$ such that
\begin{equation}
\label{eq:f}
Tf(\xi)=\chi_E(\xi)\widehat{f}\bigl(\ffi(\xi)\bigr)
\end{equation}
which completes the proof.
\end{proof}

\begin{remark}
If $T$ extends to a unitary operator from $L^2(\R^d)$ onto $L^2(\R^d)$ then $E=\R^d$
and $\ffi\,:G\to G$ is bijection and is measure preserving {\it i.e.} $|\ffi^{-1}(E)|=|E|$
for every set $E\subset G$ of finite measure. This last fact is a corollary of \cite{Si} ({\it see also} \cite{No}).
\end{remark}

Note that in this theorem, we have only used the $L^1-L^\infty$ duality to show that the operator
is a kernel operator. A slightly more evolved theorem allows to obtain this directly. More precisely,
this is a consequence of the following theorem that dates back at least to Gelfand \cite{Ge}
and Kantorovich-Vullich \cite{KV} ({\it see} also \cite[Theorem 2.2.5]{DP} or \cite[Theorem 1.3]{AT}):

\begin{theorem}
Let $(\Omega_1,\mu_1)$ and $(\Omega_2,\mu_2)$ be two $\sigma$-finite measure spaces. There is a one-to-one correspondence
between bounded linear operators $T\,:L^1(\Omega_1)\to L^\infty(\Omega_2)$ and kernels $k\in L^\infty(\Omega_1\times\Omega_2)$. This correspondence is given by $T=T_k$ where $T_k$ is defined by
$$
T_kf(\omega)=\int_{\Omega_1}k(\zeta,\omega)f(\zeta)\,\mbox{d}\mu_1(\zeta),\quad f\in L^1(\Omega_1).
$$
\end{theorem}

It follows that this proposition than essentially reduces to the results in \cite{Lu1,Lu2},
However, a non explicit condition in those papers is that $k$ be defined everywhere as it is applied to Dirac masses.

\section{The twisted convolution}

In this section, we consider the case of the twisted convolution (for background on this transform we refer to \cite{Fo}). Recall that this is defined for
$f,g\in L^1(\R^{2d})$ by
$$
f\cvh g(x,y)=
\int_{\R^{d}}\int_{\R^{d}}f(x-s,y-t)g(s,t)e^{i\pi(\scal{x,t}-\scal{y,s})}\,\mbox{d}s\,\mbox{d}t.
$$
This defines a new $L^1(\R^{2d})$ function.
Note also that this operation is non commutative.

Next, for $p,q\in\R^d$, let us define the following operator that acts on functions on $\R^d$:
$$
\rho(p,q)\ffi(x)=e^{2i\pi \scal{q,x}+i\pi\scal{p,q}}\ffi(x+p).
$$
For $f\in L^1(\R^d\times\R^d)$ we define the (bounded linear) operator on $L^1(\R^d)\to L^1(\R^d)$
\begin{eqnarray*}
\rho(f)\ffi(x)&=&\int_{\R^d}\int_{\R^d}f(p,q)\rho(p,q)\ffi(x)\,\mbox{d}p\,\mbox{d}q\\
&=&\int_{\R^d}K_f(x,y)\ffi(y)\,\mbox{d}y
\end{eqnarray*}
where $\dst K_f(x,y)=\int_{\R^d}f(y-x,q)e^{i\pi\scal{q,x+y}}\,\mbox{d}q=\ff_2^{-1}[f]\left(y-x,\frac{x+y}{2}\right)$
and $\ff_2$ stands for the Fourier transform in the second variable.

One then checks through a cumbersome computation that $\rho(f\cvh g)=\rho(f)\rho(g)$ (here the product stands for composition of operators) or, for the kernels
$$
K_{f\cvh g}(x,y)=\int_{\R^d}K_f(x,z)K_g(z,y)\,\mbox{d}z.
$$

\medskip

\noindent{\bf Question.} {\sl To what extend does this characterize the transform $f\to\rho(f)$.}
%
%
%
%
%

\end{document}